\documentclass{article}                    
\usepackage{graphicx}
\usepackage{amsmath,amssymb}
\usepackage{amsthm}

\newtheorem{remark}{Remark}
\newtheorem{theorem}{Theorem}
\newtheorem{lemma}{Lemma}
\newtheorem{corollary}{Corollary}
\newtheorem{proposition}{Proposition}
\begin{document}

\title{A Note on Lerner Index, Cross-Elasticity and Revenue Optimization Invariants %\thanks{Grants or other notes
}


\author{Alexander Kushkuley (UTD, kushkuley@gmail.com)        \and
        Su-Ming Wu (Oracle) %etc.
}

\date{November 2013}

\maketitle

\begin{abstract}
We study common properties of retail pricing models within a general framework of calculus of variations. In particular, we observe that for any demand model, optimal de-seasoned revenue rate divided by price elasticity is time invariant. We also obtain a generalization of a well known inverse relationship between price elasticity of demand and Lerner index. These invariance results are illustrated by two contrasting examples of markdown optimization and optimal continuous replenishment
%\keywords{variational methods\and (cross-)elasticity of demand \and inventory effect \and Lerner index \and markdown optimization  }
\end{abstract}

\section{Introduction}
\label{intro}
 
In order to maximize his profit a retailer will try to exercise an optimal price changing policy and an optimal inventory replenishment policy. 
We thus have a variational optimization problem
%, where a goal function is represented by a demand model 
for an  expected profit  functional defined on a 
control space of time dependent price and inventory vector-functions and 
we would like to find out what (if any) price/inventory policy guidelines can be obtained  from a general variational
formulation of multi-item revenue/profit maximization problems for arbitrary continuous demand models. In some sense this paper is an  attempt at meta-analysis of  variational 
demand models. The questions that we are posing are:
\begin{itemize}
\item[.] What demand models admit a reasonable first order optimality conditions?  
\item[.] What restrictions on price elasticity and inventory effect  are imposed by optimality? 
\item[-] What are time invariants of optimal solutions?  
\end{itemize}

Remarkably, quite a few useful insights can be obtained even in a most general setting. For example,
we observe the relation of optimality conditions to Lerner index (cf. \cite{Thijs}) and 
obtain general markdown profit/revenue invariance conditions that were established in one-item case for exponential elasticity model in \cite{Smith2010}, \cite{Smith1998},  and for constant elasticity model in \cite{Vak})
As an application we find a realistically looking closed form solution for 
constant elasticity  multivariate markdown optimization  problem with inventory effect (for the "univariate" one-item case, cf. \cite{Smith2010}, \cite{Smith1998} and\cite{Vak}) 
as well as a closed form solution for the opposite continuous replenishment case.

The authors are aware of two cases of practical applications of  markdown invariants: 
\begin{itemize}
\item[(i)] %it is shown inthat 
an optimal markdown policy for exponential demand model is to keep seasonally adjusted rate of sales constant(\cite{Smith1998}, \cite{Smith2010})
\item[(ii)] an optimal markdown policy for constant elasticity demand model  is
to keep de-seasoned revenue rate constant (cf. \cite{Vak}) 
\end{itemize} 
and probably the most important result of this study is that these
 seemingly contradictory statements are in fact 
 manifistations of a more genera invariance principle  (cf. Corollary 2 below).  

\section{Preliminaries}
\label{intro}

 A continuous profit/revenue optimization problem can be  parametrized by  the following time dependent vector variables indexed by the number of distinct items (products) (cf. e.g. \cite{Talluri}) 
\begin{itemize}
\item Inventory $ I \; \equiv \; I(t) \; \equiv \; (I_1, \cdots , I_n ) , I_i(t)  \geq 0 \; , \; i = 1, \cdots, n   $    
\item Price $ p \; \equiv \; p(t) \; \equiv \; (p_1, \cdots , p_n ),  \;   p_i(t)  \geq 0 \; , \;  i = 1, \cdots, n $ 
\item Rate of sales $ S \; \equiv \; S(t) \; \equiv \; (S_1, \cdots , S_n ), \; ,\; S_i(t) \geq 0 ,  \; i = 1, \cdots, n$   
\item Revenue rate $ R  \equiv  R(t)  \equiv  (R_1, \cdots , R_n ), \;  R_i(t)  = p_i(t) S_i(t),  \; i = 1, \cdots, n$   
\item Replenishment  rate $ \rho \; \equiv \; \rho(t) \; \equiv \; (\rho_1, \cdots , \rho_n ), \; \rho_i(t) \geq 0 \; ,      \; i = 1, \cdots, n    $  
\item Retailers cost of one item  $ c   \; \equiv (c_1, \cdots , c_n ), \; c_i \geq 0 \;     \; i = 1, \cdots, n    $   
\item Profit  rate $ P  \equiv  P(t)  \equiv  (P_1, \cdots , P_n ),   P_i(t)  = (p_i - c_i ) S_i(t),  \; i = 1, \cdots, n$
\item Lerner index   (cf. \cite{Ler})
 $ l\equiv l(t)  \equiv (l_1, \cdots , l_n ), \; l_i(t)  = \frac{p_i(t) - c_i }{p_i(t) }      \; , \;i = 1, \cdots, n    $  
\end{itemize}  
In what follows we use a component-wise multiplication of vectors omitting   indices of vector variables.
For example, the relation between profit vector $P,$ revenue vector $R$ and the vector of Lerner indices $ l $ can be written as 
\begin{eqnarray}
 P = l R 
\end{eqnarray} 
In this notation for a matrix $ M $  applied to a coordinate-wise product of two vectors $ A $ and $ B $ we have   
\begin{eqnarray}
M(AB) = (M diag( A )) B = (M diag(B)) A = M(BA)  \nonumber 
\end{eqnarray}   
where $ diag( A ) $ is a diagonal matrix with coordinates of a vector. We will
omit brackets and diagonalization signs in similar formulas.  The notation $MA$ is thus ambiguous, since it can denote a matrix $ M diag( A ) $ or a vector $ M(A) .$  It is,  however, usually clear from the context whether the result is a matrix or a vector. 
 
It is assumed that there are no "stock outs" for any of the items involved and hence that a rate of sales of 
any item is the same as a (consumer) demand for the item. 
 A standard econometric model (cf. \cite{Talluri}, \cite{Smith1998}, \cite{Smith2010}, \cite{Vak}  for one item formulation)  
postulates that the rate of demand (sales) and hence revenue and profit rates  are functions of inventory and price
\begin{eqnarray}
S = S( I, p)  \; \Leftrightarrow \; S_i(t) = S_i( I(t), p(t) ) \;  \; , \; i = 1, \cdots ,  n   \nonumber \\
R = R(I,p) \equiv Sp, \; \; P = P(I,p) \; \equiv \; S(p-c)\;  \nonumber
%\Leftrightarrow \; R_i(t) = R_i( I(t), p(t) ) \; \equiv S_i(t)(p_i(t) - c_i) \nonumber 
\end{eqnarray}
Below we will often write $ I, p $ instead of $ I(t), p(t) $ and, for example,  $ (p-c)S $ instead of $ (p(t) - c )S(p(t) ,I(t) )$ and so on. 
A total revenue or profit (expectation) is given  as an  integral of revenue (or profit) rate along some predefined time  interval $ [ 0, T ] .$ Such an integral is usually taken with respect to a beforehand estimated measure  $\sigma(t) dt$ where $ \sigma \equiv \sigma(t) $ could be a known density of some common to all items random or deterministic shock, for example seasonality effect  (see\cite{Smith1998},  \cite{Smith2010}, \cite{Vak}). Regardless of its nature,  we will call
the density $ \sigma $ "seasonality" assuming without a loss of generality that it is normalized as a probability density.
\begin{eqnarray}
\int_{0}^{T}\sigma(t) dt  \; = \; 1  \nonumber 
\end{eqnarray}

Hence,  profit expectation functional , for example,  can be written as
\begin{eqnarray}
 \int_{0}^{T} <1_n,  S(p-c )>  \sigma (t) \; \equiv \;   \int_{0}^{T} <S,  (p-c )>\sigma dt 
\end{eqnarray} 
where $ 1_n$ denotes an $n-$vector whose all coordinates are equal to $1$ and angle brackets denote the standard scalar product. 

%where 
%$\sigma(t) dt$ is some probability measure, for example $ \sigma $ could be a known %probability density of some common to all items random shock (effect)   
% (e.g seasonality, see\cite{Smith1998},  \cite{Smith2010}, \cite{Vak}).
The main subject of this study is the following constrained variational problem (cf. \cite{abs})
%\begin{eqnarray}
\begin{align}
&\texttt{maximize}   \nonumber \\
& \int_{0}^{T}   [\; \; - <c,\rho^2(t)> +  <p(t), \; S_i(I(t),p(t))>  \; ] \sigma dt  \\
&\texttt{subject to an inventory flow constraint}   \nonumber \\ 
&\frac{ dI } {dt} \; = - \; S(I(t),p(t), t) \sigma + \rho^2(t) \sigma   \\
&\texttt{and additional boundary constraints, e.g.}  \; I(0) = I_0, \; I(T) = I_T \nonumber
\end{align}
%\end{eqnarray}
 % \; p(0) = p_0, \cdots ,$ etc. 
Note that mnemonics $ \rho^2(t) $ is used to denote non-negative de-seasoned product replenishment rate.
  %or equivalently 
%\begin{eqnarray}
% \sum_{i} \int_{0}^{T}   -c_i(   \frac{ dI_i } {dt} \; + \; S_i  \sigma(t) )  +  p_i S_i    \sigma(t) dt  \\
%\frac{ dI_i } {dt} \; + \; S_i\sigma \geq 0 \; , \: i = 1, \cdots , n  
%\end{eqnarray}
For simplicity,  we assume fixed  cost of replenishment  $ c .$ The density  $ \sigma(t) $ is assumed to be known in advance and the sought for "control" functions are replenishment $ \rho \; \equiv \rho(t) $, price
 $p \; \equiv p(t)$    
and inventory $ I \equiv I(t)  .$ 
When the items are not replenished  
($ \rho \equiv  0 $),  the problem is reduced to the so called markdown optimization problem (cf. \cite{Smith1998}) 
\begin{align}
&\texttt{maximize} \nonumber \\
&\int_{0}^{T}    <p(t), \; S_i(I(t),p(t))> \sigma \;  dt  \\
&\texttt{subject to inventory constraint} \nonumber \\ 
&\frac{ dI } {dt} \; = - \; S(I(t),p(t), t) \sigma    \\
& \texttt{and markdown  boundary constraints}  I(0) = I_0, \; I(T) = 0 \nonumber  
\end{align}
\begin{remark}
In this mostly qualitative study we consider only simplest boundary constraints.
We also ignore various important business rules that usually impose relations between item variables
e.g.  $ p_1 + p_2 < p_3 $ or $ I_5 < 100 $. See, however, concluding remarks in section 7.
\end{remark}

% One can also consider also a case of "unconstrained" revenue optimization (5)  ignoring any cots and inventory constraints ( $c(t) \equiv 0 $). 
%\begin{eqnarray}
%max \int_{0}^{T}    <p(t), \; S(I(t),p(t))> \sigma \;  dt  
%\end{eqnarray}
\section{Cross-elasticity and cross-inventory effects}
\subsection{Definitions}
It is customary to express effect of price changes on demand in terms of price elasticity (cf. \cite{Talluri}).     
Recall that cross-elasticity matrix $\Gamma$ of our item bundle  is defined as a Jacobian matrix of $ S $ with respect to $ p $
multiplied by the diagonal matrix $ S^{-1} $ on the left and diagonal matrix $ p $ on the right (cf. e.g. \cite{Thijs})
\begin{eqnarray}
  \Gamma =   \parallel S_i^{-1} \partial {S_i}  / \partial{p_j} \;  p_j \parallel \; \equiv \; 
  diag( S^{-1} ) \frac{ \partial S }
  { \partial p }  diag ( p )        
\end{eqnarray}
By analogy we introduce also a cross-inventory effect matrix 
\begin{eqnarray}
  \alpha =   \parallel S_i^{-1} \frac{ \partial {S_i} }  { \partial{I_j} } \;  I_j \parallel \; \equiv \; 
  diag( S^{-1} ) \frac{ \partial S } { \partial I }  diag ( I )        
\end{eqnarray}
Omitting  diagonalization signs, we have, in other words 
\begin{eqnarray}
\frac{\partial S }{\partial p }  = S \Gamma p^{-1}  \\
 \frac{ \partial S } { \partial I }  = S \alpha I ^{-1}    
\end{eqnarray}

We recall also (see e.g. \cite{Talluri}, \cite{Thijs}) that elasticity effect $ \gamma_{i,j} $  measures change in relative demand of an item $ i $ due to the relative price change of an item $ j  .$  Since self-effect $ \gamma_{i,i} $ probably affects demand for the item more than demand for other items,   it is reasonable to assume that $ \Gamma $ is diagonally dominant, i.e.
$ \gamma_{i,i} < 0 $  and $ \sum_{j} |\gamma_{i,j} | < |\gamma_{i,i} | ).$ 
(cf. e.g. \cite{Thijs}, \cite{horn}). 

\paragraph{}
Quite similarly,   a cross-inventory effect $ \alpha_{i,j} $ is a change in relative demand for item $ i $ caused by relative change in inventory of item $ j $  and we  mention in passing that  inventory effect can be viewed as a  formal generalization of a notion of "demand transference".    
 When an item $ i $ becomes scares (or is removed from the shelves) a demand for the item is  transferred to other items
$ j_1, j_2, \cdots .$ The amount of this "transference" can be  quantified by the $i , j_k $ element of the matrix $ \alpha .$
To be more precise we can define the demand transference of item $i$ to
item $j$ as 
\begin{equation}
d_{i,j} = \lim_{I_i \rightarrow 0}  \left(  S_i^{-1} \frac{ \partial {S_i} }  { \partial{I_j} } \;  I_j  \right)  \; \equiv \; \lim_{I_i \rightarrow 0}  \alpha_{i,j}
\end{equation}  

\subsection{Constant and Exponential elasticity}

As a standard example of the defining relations (9-10) we mention a log linear model
\begin{equation}
\log S = \log S_0 + \alpha \log I + \Gamma \log p 
\end{equation}
or equivalently
\begin{equation}
\log R = \log S_0 +  \alpha  \log I  + ( 1 + \Gamma)  \log p 
\end{equation}
When all the matrix coefficients of $ \alpha $ and  $ \Gamma $ do not depend on time this is a so called constant elasticity model or SCAN*PRO with inventory effect (the "pure" SCAN*PRO 
requires $ \alpha $ to be zero cf. \cite{Talluri}).
For example (cf. \cite{Vak}) one-dimensional SCAN*PRO with inventory effect can be specified as
\begin{equation}
S = S_0 I^{\alpha} p ^{\gamma} 
\end{equation}
 for some constant base demand rate $S_0$ and "effect numbers"  $\alpha$ and $\gamma.$
\paragraph{}
   By analogy,  an exponential demand model with constant coefficients (cf. e.g. \cite{Smith1998}, \cite{Smith2010}) is specified as  
\begin{equation}
\log S \; = \log S_0\; + \; \alpha \log I + \Gamma p 
\end{equation}
for some constant matrices  $ \Gamma $ and  $\alpha, $ and in one dimensional case as
\begin{eqnarray}
S \; = \; S_0 I^{\alpha} e^{p\gamma} 
\end{eqnarray}

As it should be expected, the exponential model (15) is not constant elasticity model. 
It is, however,  easy to see that  elasticity matrix  $  \tilde{\Gamma} $ for exponential model (15)  is equal to   
\begin{equation}
 \tilde{\Gamma}   = \Gamma p \; \equiv \; \Gamma diag( p ) 
\end{equation}
For example in one dimensional case the elasticity of exponential model with constant coefficients is  
\begin{eqnarray}
\tilde{\gamma} \; = \; S^{-1} \frac {\partial S}{\partial p } p = \gamma p 
\end{eqnarray}

\section{First Order Conditions}
It is not hard to write down necessary optimality conditions for the problem (3-4) (in a different context similar computation was done 
in \cite{Thijs}). 
 
First, note that the optimization problem (3-4) is equivalent to 
\begin{eqnarray}
max \int_{0}^{T}   [\; \; - <c,\frac{ dI } {dt} \; + \; S(I(t),p(t), t)>  +  <p(t), \; S_i(I(t),p(t))> \; ]  \sigma dt   \nonumber \\
\frac{ dI } {dt} \; + \; S(I(t),p(t), t) \sigma \;  \geq \;0 \; \; \; 
\end{eqnarray}
so that the vector function $ \rho^2 $ plays a role of a dummy variable.
Therefore (see e.g. \cite{deim}, chapter 9) under some general conditions on functions involved a necessary optimality condition for the problem (19) can be written down in terms of "Lagrange multiplier" functions $ \lambda_i \; \equiv \; \lambda_i(t) $  such that on the optimum trajectory ( I(t), p(t) ), if it exists, either  $ \lambda_i = 0 $ or  $ \frac{ dI_i } {dt} \; + \; S_i(I(t),p(t), t) \sigma = 0 $ and $ \lambda_i <= 0 $ for every item $i$. Using retail business terminology, one can say that this "Lagrange multiplier"  alternative distinguishes between continuous replenishment and markdown (see below).  It is more convenient for our purposes, however, to stick to the original formulation (3-4). Introducing Lagrangian multiplier vector function $\lambda$ ( cf. e.g \cite{gelf} ) we can rewrite (3,4) as
\begin{eqnarray}
\int_{0}^{T} L dt \equiv  \int_{0}^{T}  [-<c, \rho^2 > \sigma + <p,S>  \sigma + <\lambda,  \frac{dI}{dt}  + S\sigma - \rho^2 \sigma> ]   dt 
\end{eqnarray}  
and hence we have 
\begin{theorem}
The first order optimality (equilibrium) conditions for the problem (3-4) (Euler-Lagrange equations) are 
%\begin{eqnarray}
%\dot{\lambda} = \sum_i ( - p_i + \lambda_i + c_i ) \frac{\partial S_i}{\partial I_j } \sigma \; , \; j = 1,\cdots , n  \\  
%\sum_i  \delta_{ij} S_i +  ( p_i - \lambda_i - c_i ) \frac{ \partial S_i}{\partial p_j }  \; = \; 0  \; , \; j = 1,\cdots , n 
%\end{eqnarray}
%or in matrix notations
\begin{align}
& (c+\lambda) \rho \; = \; 0  \\
& p + \lambda   = - \left[ \left( \frac{\partial S }{\partial p } \right) ^T \right] ^{-1}  S  \; \Leftrightarrow \; \lambda   =   -p R^{-1}  \left( 1 + \left( \Gamma^T   \right)^{-1}  \right) R \\
& \dot { \lambda } =  \sigma \left( \frac{\partial S }{\partial I } \right) ^T   ( p + \lambda )    \; \Leftrightarrow \; \dot { \lambda } = \sigma I^{-1} \alpha^T S ( p + \lambda ) \; \equiv \;  \sigma I^{-1}  \alpha ^T \left ( \Gamma^T   \right) ^{-1}R \\ 
& \lambda_i  = -c_i \; \; \texttt{or}   \; \;\rho_i = 0, \; \lambda_i <= 0, \; i = 1,\cdots ,n 
% \Leftrightarrow  p_i - c_i = -\left[ \left[ \left( \frac{\partial S }{\partial p } \right) ^T \right] ^{-1}  S \right]_i    \texttt{and}
%\sigma \left( \frac{\partial S }{\partial I } \right) ^T   ( p  -c ) 
\end{align}
\end{theorem}
Proof. The equations (21-23) are Euler-Lagrange equations for the Lagrangian $L$ given by (20) while (24) follows from "Lagrange multiplier" alternative (21) that was discussed above.
\paragraph{ }

In one-item case (22) boils down to $ \lambda = - ( 1 + \frac{1}{\gamma} ) p ,$ and an optimality condition alternative is reduced to 
\begin{corollary} 
Alternative optimality conditions for one-item problem (3-4)
are 
\begin{itemize} 
\item[(i)] An item is continuously replenished, $\alpha = 0 $ and price is fixed at  $ p = \frac{\gamma}{\gamma + 1} c.$ For this to make sense we need of course inequality  $ \frac{\gamma}{\gamma + 1}  \geq 1 $ to hold 
\item[(ii)] an item is being (optimally) marked down ($ \rho = 0 $) and therefore the inequality $ \frac{\gamma + 1}{\gamma} \geq 0 $ must hold
\end{itemize}
\end{corollary}
\begin{remark}
 A condition (i) above is equivalent to Lerner inverse elasticity rule $ \gamma = - p / (p-c) $  
\end{remark}  
\begin{remark} A practically useful "coincidence" is that inequalities imposed on elasticity by Corollary 1 in either case are equivalent to $ \gamma \leq -1 $ which is a  usual assumption for an elastic item (cf. e.g. \cite{Talluri}, \cite{Vak}).   We will discuss a multi-item version of Corollary 1 further on. 
\end{remark}
\begin{remark} The corollary implies that optimal continuous replenishment  equilibrium of one item cannot exist if inventory effect $\alpha$ is non-zero. An intuitive meaning of this condition is that when an item can be replenished "at will" the inventory effect should not matter. As we will see below, this is not necessarily true  in multi-item case, since it might not be possible to eliminate cross-inventory effects by 
replenishment. 
\end{remark}
%\begin{eqnarray}
%\rho_i = 0 \; \; \texttt{or} \; \;    p_i - c_i = -\left[ \left[ \left( \frac{\partial S }{\partial p } \right) ^T \right] ^{-1}  S \right]_i                     
%\end{eqnarray}
It  follows from the standard calculus of variations (cf. e.g. \cite{gelf}) that  Hamiltonian associated to the Lagrangian function $ L $ in (20) is
\begin{eqnarray}
H \; = \; < \frac{ \partial L } { \partial \dot{I} },  \dot{I} > \; - \; L  \; = \;  <\lambda, \frac{dI}{dt} > 
+ <c,\rho^2> \sigma -   <p,S>  \sigma   \nonumber \\
-<\lambda,  \frac{dI}{dt}  + S\sigma - \rho^2 \sigma> \nonumber \\ 
 = \; <c+\lambda, \rho^2> \sigma  \; - \; <p + \lambda,S>\sigma \; \; \; \; \;\;\;\;\;\; \nonumber
 \end{eqnarray}  
Therefore, the quantity 
\begin{equation}
<c+\lambda, \rho^2>  - <p + \lambda,S> 
\end{equation}
is an optimal sales invariant in a sense that it does not change on optimal price/inventory trajectory. We will discuss some  practical implications of this fact below. It is, however, somewhat satisfying to verify time invariance of the expression (25) directly.  Computing
time derivative of (25) and applying (23) we get
\begin{eqnarray}
\frac{d}{dt}  ( <c+\lambda, \rho^2> - <p +\lambda, S > )  =  \nonumber \\
\sigma < \left( \frac{\partial S }{\partial I } \right) ^T   ( p + \lambda ), \rho^2> +  
2<c+\lambda, \rho \dot{\rho} >   \nonumber \\ 
-    < \frac{dp}{dt}  + \sigma \left( \frac{\partial S }{\partial I } \right) ^T   ( p + \lambda )  , S > - <p + \lambda, 
\frac{\partial S}{\partial I } \frac{dI}{dt}  + \frac {\partial S}{\partial p} \frac{dp}{dt} > 
\nonumber \\
\; = \;   \sigma < \left( \frac{\partial S }{\partial I } \right) ^T   ( p + \lambda ), \rho^2> +  
2<(c+\lambda) \rho,  \dot{\rho} >   \nonumber \\
- <\frac{dp}{dt}, S > -  \sigma <\left( \frac{\partial S }{\partial I } \right) ^T   ( p + \lambda )  , S >  \nonumber \\  -   <  \left( \frac{\partial S }{\partial I } \right) ^T   ( p + \lambda )  , \frac{dI}{dt} > -
<  \left( \frac {\partial S}{\partial p} \right) ^T ( p + \lambda) , \frac{dp}{dt} >  
\end{eqnarray} 
Replacing here $ dI/dt $ by the right hand side of (4) and applying (21) we observe that (26) reduces to (22) demonstrating once gain that (25) is a constant (time invariant) on optimal trajectory. Now, substituting the expression for Lagrange multipliers    $ \lambda $  given by (22) into our invariant (25) and using definition of elasticity matrix (9)  we obtain
\begin{eqnarray}
<c+\lambda, \rho^2>  - <p + \lambda,S>  \nonumber \\
= \;  < c - p  - S^{-1} \left( \Gamma^T \right) ^{-1} p  S, \rho^2 >  + < S^{-1} \left( \Gamma^T \right) ^{-1} pS, S > \nonumber \\
=  \; < (c - p )S - \left( \Gamma^T \right) ^{-1} p  S, \frac{\rho^2}{S} > + <\left( \Gamma^T \right) ^{-1} p  S, 1_n > \nonumber \\
= \; - < p - c,  \rho^2  > + < R,  \Gamma^{-1}  ( 1_n -   \frac{\rho^2}{S}  ) >
\end{eqnarray}
We have established, therefore,
\begin{theorem}
The quantity (27) is time invariant on  optimal trajectory of the problem (3-4) if such an optimal trajectory exists.
\end{theorem}
\section{Continuous replenishment versus markdown }
In practical terms, Theorem 1 implies that optimal price/inventory control requires every item to be either continuously replenished or marked down. It is reasonable to assume that a bundle of related items is marked-down simultaneously. 
Therefore, in what follows we consider two opposite cases 
\begin{itemize}
\item[(CR)] optimal continuous replenishment -  the replenishment rate is non-zero for every item 
( $ \rho_i \neq 0 , i = 1, 2, \cdots , n ) ,$ 
\item[(MD)]  mark down optimization (MDO) - all items  are marked down ( $ \rho_i = 0  , \; i = 1, 2, \cdots , n $ )
\end{itemize}
\begin{remark}
A continuous replenishment case (CR) roughly corresponds to what is known in practical applications as a regular price optimization (RPO). The anecdotal evidence in such cases is that regular (de-seasoned, non-promotional)  prices remain  unchanged for relatively long periods of time. This is especially true for staple or grocery items.    
\end{remark}

Both cases have some common features 
\begin{corollary}
\
%\paragraph{}
\begin{itemize}
\item[(i)] The quantity  $- < R,  \Gamma^{-1} 1_n > $ is invariant on optimal price/inventory trajectory in both cases (CR) and (MD)
\item[(ii)] In case (CR)  this quantity is
equal to a total profit rate, i. e.  $ < p - c , S > \; =  \;  - < R,  \Gamma^{-1} 1_n > $
\end{itemize}
\end{corollary}
     
Proof. In case (MD) the first term of invariant (25) vanishes because $ \rho = 0 .$
In case (CR) the same term vanishes because $ \lambda = - c .$ For the same reason, in case (CR) the invariant (25) is equal to $- <p-c, S> .$    Finally it follows from (27) that the second term of (25) is always equal to $ < R,  \Gamma^{-1} 1_n > .$
\paragraph{}
Recalling the definition of exponential model given above in section 4 we get 
\begin{eqnarray}
< R,  \tilde{\Gamma}^{-1} 1_n > = <R, p^{-1}\Gamma^{-1} 1_n> = 
<S, \Gamma^{-1} 1_n> \nonumber
\end{eqnarray}
and hence, the following
%The authors are aware of two cases of practical applications of  markdown invariants 
%\begin{itemize}
%\item[(i)] It is shown in \cite{Smith} that the optimal markdown policy for one dimensional exponential demand model is to keep seasonally adjusted rate of sales constant. In our notation that can be written down as $ S = const $  
%\item[(ii)] On the other hand, the optimal markdown policy for one dimensional constant elasticity demand model  is
%to keep de-seasoned revenue constant, i.e. $ R = const .$ (cf. \cite{Vak}) 
%\end{itemize} 
\begin{corollary} 
\
\begin{itemize} 
\item[(i)] For exponential model (15) with constant exponent matrix $\Gamma$ the quantity
$ \;< S, \Gamma^{-1}1_n> $ is constant on optimal price/inventory trajectory in both cases (MD) and (CR)
\item[(ii)] In particular, for one dimensional exponential model an optimal price/inventory policy is to maintain deseasoned sales constant  (cf, \cite{Smith1998}, \cite{Smith2010})
\item[(iii)] For one dimensional constant elasticity model an  optimal price inventory policy is 
to maintain deseasoned revenue constant (cf. \cite{Vak} )
\end{itemize}
\end{corollary}
\paragraph{}
From the general optimality conditions (22-24) we also infer a generalization of one-dimensional  Corollary 1 (see also \cite{Thijs} )
\begin{corollary}  
\
%\begin{itemize}
%\item[(CR Case)] 
In CR case, necessary optimal equilibrium conditions for a demand model with continuous replenishment are as follows
\begin{itemize} 
 \item[(i)] At any time  on optimal price/inventory trajectory the revenue vector  is an eigenvector of 
 $ \Gamma^T l $ with eigenvalue $-1.$ In other words,  a generalized  Lerner inverse elasticity rule holds    
\begin{align}
\Gamma^T l R = - R \; \Leftrightarrow \;  \Gamma^T P = -R  \; \Leftrightarrow   \;    
   l = - R^{-1} (\Gamma^T)^{-1} R 
\end{align} 
Clearly, for this to make sense it is necessary that 
\begin{equation}
\left( \Gamma^T \right) ^{-1} R \geq -R   \; \Leftrightarrow \; \left( 1 + \left ( \Gamma^T \right) ^{-1}\right) R \geq 0
\end{equation} 

\item[(ii)] The inventory effect matrix $ \alpha $ must satisfy the following degeneracy condition
\begin{equation}
\alpha^T S (p-c)  = 0 \; \Leftrightarrow \; \alpha^T P =   0  \; \Leftrightarrow \;  \alpha^T l R = 0
\end{equation}
\end{itemize} 
Necessary conditions for mark down optimality (case MD) are equations (22-23) and inequality (29).
% $ \left( 1 + \left ( \Gamma^T \right) ^{-1}\right) R \geq 0 $ 

%\end{itemize}

\end{corollary}
Proof. The formulas (28-29) directly follow from (22-24) and definitions of inventory effect $ \alpha $ and elasticity $\Gamma $ (9-10).

\begin{remark}
The condition (29) is a matrix equivalent of the negative elasticity condition (ii) of Corollary 1 that was discussed in Remark 3. This more general condition seems to be as reasonable as its one-dimensional prototype in a sense that matrix with large in absolute value negative diagonal behaves like a negative number at least when applied to a profit vector $ R .$  We will call such a matrix "highly negative" and will assume below that elasticity matrix $ \Gamma $ satisfies this condition. Note that diagonal dominance is in general a weaker condition 
\end{remark}
\begin{remark}
The condition (30) is a matrix equivalent of the condition (23) of Theorem 1  that was mentioned in Corollary 1 and Remark 4. It implies  that in one-dimensional (one-item) CR case the demand $ S $ does not depend on inventory
at all,  $ \alpha \equiv  0 .$ In general, the condition (30) means that in CR case the first order optimality conditions cannot be satisfied unless  
\begin{eqnarray}
\forall i : \sum_{j=1}^{n} \Delta^{i}_{j} S_j (p_j - c_j) = 0 
\end{eqnarray}
where  $ \Delta^{i}_{j}S_j$ denotes a change in demand of an item $ j $ due to change in inventory of an item $ i .$
This seems to be intuitively clear - if unlimited replenishment is allowed then at an equilibrium, a gain (loss) from change in inventory of an item will be compensated for by corresponding loss (gain)  in other items.
\end{remark}

\section{Closed form solution examples} In this section we will demonstrate that Theorem 2 and Corollaries 3,4 can be, at least in principle, used to find realistic optimal price inventory policies.  In what follows we assume that both inventory effect matrix and elasticity matrix are constant. Cases CR and MD will be considered separately. However in either case we will use the same heuristics that we will briefly describe now.
\paragraph{•}
According to Corollary 2 (i) the quantity
\begin{align}
\sum _i a_i R_i  \nonumber
\end{align}
 must be constant on optimal trajectory. Here $ a_i $ is a sum of elements of the $i-$th row of $ \Gamma^{-1}  , \; i = 1, \cdots, n .$ 
In case of constant elasticity, all the numbers $a_i $ are themselves constant and therefore a price/inventory trajectory that keeps all individual item revenues $ R_i $ constant has a chance to be optimal. This heuristics is quite restrictive and precisely because of that it  "selects" solutions that can be investigated  by hand.  We will impose additional conditions on effect matrices $ \alpha $ and $ \Gamma $ as we go -  after all, the purpose of this exercise is  to develop intuition in problems of this kind and to show that conservation principle established in Corollary 2 
  can serve as a guideline in revenue optimization.

\subsection{Optimal solution: continuous replenishment case} 
An optimal pricing  recipe for one dimensional case is essentially supplied by Corollary 1. Assuming that there is no inventory effect ( $\alpha = 0 $ ) maintain constant price $ p  =  \frac{\gamma} { \gamma + 1 } c $ and keep inventory constant ( $ dI/dt = 0 $ ) by full replenishment  $ \rho^2 = S$ 
As simple as it is, this solution is in a good agreement with commonly used polices for staple items. A similar approach can be taken in a multi-item case. 
We will consider for simplicity a constant effect model (13).  As was explained above  the condition 
(ii) of Corollary 2 will be satisfied if all the components of the profit vector are kept constant. For constant effect model we thus have 
\begin{align}
& P = S( p - c ) = const   
\end{align}  
Applying (28) and (1)  we turn (32) into a system of equations
\begin{align}
& \alpha \log I  + ( 1 + \Gamma )\log p    = \; \log R   \\
& (\Gamma^T)^{-1}R \; = \; -lR \; \equiv \; -P
\end{align}
%and (32) can be rewritten as 
%\begin{eqnarray}
%\alpha \log I + ( 1 + \Gamma ) \log p = \log ( - (\Gamma^T)^{-1} P ) 
%\end{eqnarray}
Since $ P, \alpha, \Gamma $ are all constant, it follows from  (34) that $ R $ and 
$ l $ are constant. Since price $ p $ is completely determined by Lerner index  $ l $ it is also constant and then from (33) it follows that 
inventory $ I $ is constant as well. Hence we assume that $ I(0) = I(T) $ and $ \rho^2 = S .$   

Let's see if we can find a solution of (3-4) that satisfies (33). 
%First of all, the condition (28)  implies that  price must be also constant since  it is completely determined by the  Lerner index, $ l = ( p - c )  / p  \; \Leftrightarrow \; p = c /( 1-l ) .$   As in one dimensional case we can also keep inventory constant by fully replenishing every item involved, i. .e by setting $ \rho^2 = -S .$ That will satisfy the inventory flow equation (4). 
To satisfy the optimality condition (30) we will require  inventory effect matrix to be degenerate of maximal rank n - 1 and we will assume that there is a positive  vector $ P  $  such that $ \alpha^T P  = 0  $ so that P is determined up to a positive scalar multiplier.  Then assuming that $ \Gamma $ is highly negative matrix,  $\Gamma^T P   = -R $  for some positive vector $ R  $ and hence there is positive vector $ l $ such that $ P = l R .$  Setting $ p = c / ( 1 - l ) $  
we need to find $ I $ from the equation 
\begin{align}
& \alpha \log I  = -  ( 1 + \Gamma )\log ( c / ( 1- l ) )    +   \log R  
\end{align} 
The problem is, that matrix $ \alpha $ is degenerate. However, vector $ R $ is defined up to a positive scalar multiplier $ r $ and in fact, we  have to find positive vector $ I $ such that 

 \begin{align}
& \alpha \log I  =  -  ( 1 + \Gamma )\log ( c / ( 1- l ) )    +   \log R  + (\log r) 1_n
\end{align}   
where positive number $ r  $ can be appropriately chosen. Therefore, assuming in addition that vector $ 1_n $ does not belong to the range of $ \alpha $   we will make sure that equation (36) has a solution. We have proved the following 
\begin{proposition}
Suppose that constant inventory effect matrix $ \alpha $ is degenerate of maximal rank $ n - 1 $ and that 
the vector $ 1_n  $ does not belong to the range of $ \alpha .$  Suppose in a adition that $ \alpha^T P = 0 $ for some positive vector $ P $ and that constant elasticity matrix $ \Gamma $ satisfies conditions of Remark 6. Then if  boundary conditions  allow to maintain constant inventory ($ I(0) = I(T) $) the  optimization problem (3-4)  has a constant price, constant inventory solution that satisfies necessary optimality conditions with revenue, price and inventory  satisfying (determined by) (33-34). 
\end{proposition}

%We will also assume that there exists a positive diagonal matrix $ \delta $ such that 
%\begin{equation}
%\alpha   = \delta \alpha^T \delta^{-1} , \;  \delta > 0  
%\end{equation} 
%\begin{remark}
%the condition (34) is not too restrictive. For example it is satisfied when $ \alpha $ is symmetric or
%when $\alpha $ is a two by two matrix with off-diagonal elements of the same sign.
%\end{remark}
%Using this condition on $ \alpha $ we can set our constant  inventory to satisfy equation $ \log I = \delta_1 P $ for some scalar multiple $ \delta_1 $  of $ \delta .$ 
%This will eliminate the inventory term in (33), because $ \alpha \log I =  \delta \alpha_1^T \delta_1^{-1} \delta_1 P = 0 $ and we get the following formula for constant price 
%\begin{equation}
%\log p  \; =  \;  ( 1 + \Gamma ) ^{-1}\log R 
%\end{equation}
%In  practical computations with elasticities, prices are usually normalized to be less than 1 (divide by full price), so that left hand side of (35) is negative. Hence we must have 
%
%\begin{align}
%& - (\Gamma^T)^{-1} P \geq 1_n  \nonumber \\
%&  \alpha^T P = 0
%\end{align}
%
%
% 
%  

 \subsection{Markdown optimization}

We will further illustrate results of section 5 by directly computing  a closed form solution for a specific constant effect multivariate markdown optimization problem.  
% \begin{corollary}
% In one dimensional case, the quantity
%\begin{equation}
% R \gamma^{-1} = const
%\end{equation}
%remains invariant on optimal markdown trajectory
% \end{corollary}
% The authors are aware of two cases of practical applications of  markdown invariants 
%\begin{itemize}
%\item[(i)] It is shown in \cite{Smith1} that the optimal markdown policy for exponential demand model is to keep seasonally adjusted rate of sales constant. In our notation that can be written down as $ S = const $  
%\item[(ii)] The optimal markdown policy for constant elasticity demand model  is
%to keep de-seasoned revenue constant, i.e. $ R = const .$ (cf. \cite{Vac}) 
%\end{itemize} 
%
%These statements are not contradictory and moreover,  both of them follow from (36).  Indeed,  according to  (36), revenue rate times inverse price elasticity remains constant.  Therefore in constant elasticity case the revenue rate itself  must be constant as is claimed by  (ii). On the other hand, for exponential demand model,  the inverse elasticity according to  (25) is $ p^{-1} \gamma^{-1} $ 
%and our markdown invariant is $ ( R / p ) / \gamma \equiv S / \gamma $  as  required by  (i).
%As an example, we consider
\subsection{Constant effect one-item markdown  (cf. \cite{Vak}) }
It follows from Corollary 1 that if demand does not depend on inventory, i. e. $\alpha = 0 $
 then optimal price is constant. Assuming then that $ \alpha \neq 0 $ the Corollary 3 (ii) implies that   
  for some constant $C$
\begin{eqnarray}
I^{\alpha} p^{\gamma+1} = C  \; \Leftrightarrow \;  p = C I^{1/\theta}
\end{eqnarray}
where 
\begin{equation}
\theta = - (\gamma + 1 ) / \alpha
\end{equation}

Substituting (37) into inventory flow equation (6) we get after simple calculations 

\begin{eqnarray}
p(t) =  p_0 \left( 1 -  \hat{\sigma}(t)      \right)^{1/(1 + \theta)}  \\
 I(t) = I_0 \left( 1 -  \hat{\sigma}(t)      \right)^{\theta/(1 + \theta)}   
\end{eqnarray} 

where $ I_0,  p_0 $  are initial price and inventory,  

\begin{eqnarray}
  I(0) = I_0, \; \;  I(T) = 0   \\
  p_0 = p(0)
\end{eqnarray}
 and $\hat{\sigma}(t)  $
is a cumulative seasonality 
\begin{eqnarray}
\hat{\sigma}(t) =   \int_{0}^{t}\sigma(x) dx   \nonumber 
\end{eqnarray}  
(see \cite{Vak} for details).  Here this result is extend  to a multivariate case.

\subsection{Multivariate MDO for constant effect model}
Let $ \tau \equiv  \tau( t ) =  1 -  \hat{\sigma}(t) $ and let $a , \mu$ be a two $n$-vectors with positive coordinates. We are looking for a solution to the problem (5-6) in the following form
\begin{equation}
p = p_0 \tau^{\mu}, \; I = I_0 \tau^a
\end{equation}
We need to find conditions on $ a,b, I_0, p_0 $ that satisfy equations (22-23), (6) and the additional condition $ \lambda \leq 0 .$
Again, as was explained above, for constant elasticity, constant inventory effect model (13), the invariance condition (i) of Corollary 2  will be satisfied if we impose a 
stronger condition
\begin{equation} 
 R_i = const, \;   i = 1,2,\cdots , n  
 \end{equation} 
It is easy to see that in terms of parametrization (43) this condition is equivalent to  
\begin{equation}
\alpha a + (1 + \Gamma) \mu = 0 \; \Leftrightarrow \; \mu = - (1 + \Gamma)^{-1}\alpha a
\end{equation}    
and 
we can rewrite (6) as follows
\begin{equation}
\log \left( - \frac{dI}{dt} \right) \; = \; = \log( R/p)  + \log \sigma \; \equiv \; 
\log R - \log p + \log \sigma  
\end{equation}
Substituting (43) into (46) we get
\begin{eqnarray}
\log a + \log I_0 + (a-1) \log \tau + \log \sigma = \log R - \mu \log \tau + \log \sigma  - \log p_0  \Rightarrow \nonumber \\ 
 (a+\mu-1) \log \tau = \log R - \log a  - \log I_0 - \log p_0  \; \; 
\end{eqnarray}
and therefore 
\begin{eqnarray}
a + \mu = 1 \\
p_0 I_0 a \; \equiv \; p_0 I_0 \frac {\theta} {1 + \theta}  \; = \;  R 
\end{eqnarray}
where in accordance with (48) a positive vector $ \theta $ is chosen in such a way  that 
 $ \mu = 1 / ( 1 + \theta ) $  and $ a = \theta / (\theta + 1)  \; \equiv \; \mu \theta ;$

 \begin{remark}
 It is worth mentioning that (49) implies that only a fraction of the initial potential revenue of $ p_0 I_0 $ can be recovered by markdown. Note that because our measure $ \sigma dt $ is normalized as explained in section 2, the above expression is actually a full markdown revenue
 \end{remark}
 Using expressions for  $ a $ and $ \mu $ introduced above, we rewrite (45) as follows  
 \begin{eqnarray}
 \frac{1}{\theta} \alpha^{-1} ( 1 + \Gamma ) \mu \; = \; -\mu  \; \Leftrightarrow  \;
 ( 1 + \alpha \theta ) \mu = -  \Gamma  \mu \\
 \mu - \alpha^{-1} (1 + \Gamma)\mu = \mu + \mu \theta  \Leftrightarrow 
 [1_n - \alpha^{-1}(1 + \Gamma)]\mu = 1_n
\end{eqnarray}    
It is quite easy to see that (51) can be satisfied by a large class of diagonally dominant "highly negative" matrices $ \Gamma$ (see for example section 6.5 below). 
Next we need to look at Euler-Lagrange equations (22-23). The equation (22) can be rewritten as follows
\begin{eqnarray}
p + \lambda = - \frac{p}{R} \left( \Gamma^T \right) ^{-1} R \; \Leftrightarrow 1_n + 
\frac{\lambda}{p} \; = \; -R^{-1} \left( \Gamma^T \right)^{-1} R  
\end{eqnarray}

Since the right hand side of (52) is constant there is a constant vector $ C $
such that  
\begin{eqnarray}
\lambda \; = \;  C p  \\
1_n + C \; = \; -R^{-1} \left( \Gamma^T \right)^{-1} R   \nonumber \\
C = - R^{-1} \left( 1  +  \left ( \Gamma^T \right ) ^{-1} \right) R    
\end{eqnarray}
and therefore $ \lambda $ is less than zero as long as $ \Gamma$ is "highly negative" (cf.  Remark 6).

Now we substitute price and inventory defined by (43) into (23) and
 transform the result of this substitution by a following chain of equivalences
\begin{eqnarray}
C \dot{p}   \; = \; -\sigma I^{-1} \alpha^T (R/p)(1 + C ) p  \\
-C \mu p_0 \tau^{b - 1} \sigma \; = \; -\sigma I_0^{-1} \tau^{-a} \alpha^T R ( 1 + C ) \\
 -C \mu p_0 \; = \; -I_0^{-1}  \alpha^T \left( \Gamma^T \right) ^{-1} R  \\
-C \mu p_0 I_0 \; = \;   -\alpha^T \left( \Gamma^T \right) ^{-1} R \\
RC(1/\theta)   \;  = \;  \alpha^T \left( \Gamma^T \right) ^{-1} R \\
-\frac{1}{\theta}   R    R^{-1} \left ( 1 + \left( \Gamma^T \right)^{-1}  \right) R   \; = \;   \alpha^T \left( \Gamma^T \right) ^{-1} R \\
(\alpha^T)^{-1} \frac{1}{\theta}  \left( 1 + \Gamma^T \right) \left( \Gamma^T \right)^{-1} R  \; = \; -\left( \Gamma^T \right)^{-1}  R  \\
( \theta \alpha^T + 1 ) \left( \Gamma^T \right) ^{-1} R = -R
\end{eqnarray} 
It follows then, that  for (43) to be a solution of (22-23) and (6), the vector of revenue $ R $ must satisfy the condition (62). 
We need now a simple and well known fact from linear algebra (cf. e.g. \cite{horn} )  
\begin{lemma}
If $ A, B $ are real invertible matrices s. t. matrix $ AB $ has a real eigenvalue $ \lambda ,$
then the matrix $ A^T B^T $ has the same eigenvalue.   
\end{lemma} 
Here is a proof for completeness. We have $ ABv = \lambda v $ for some vector $ v ,$ which is the same as  $ (B - \lambda A^{-1} ) v =  0 .$ Hence, the polynomial
$ \det(  B - t A^{-1} ) \; \equiv \; \det(  B^T - t (A^T)^{-1} ) $ has a real root $ \lambda .$ Therefore $ A^T B^T w = \lambda w $ for some vector $ w .$
  
According to  (50)  the matrix  $ \frac{1}{\theta} \alpha^{-1} ( 1 + \Gamma )  $ has a positive real eigenvector $ \mu $ with eigenvalue $-1 ,$ and therefore by the above mentioned  Lemma the matrix   $ (\alpha^T)^{-1} \frac{1}{\theta}  \left( 1 + \Gamma^T \right) $ has a real eigenvector, call it $ V $  with eigenvalue -1 as well

\begin{equation}
 (\alpha^T)^{-1} \frac{1}{\theta}  \left( 1 + \Gamma^T \right) V = - V \;
 \Leftrightarrow \; \Gamma^T V \; = \; -( 1 + \theta \alpha^T ) V
\end{equation}
We will assume now that there is a positive diagonal matrix $ \delta $ such that
$ \delta \Gamma \delta^{-1} \; = \; \Gamma^T \nonumber .$
 In addition we will assume that inventory effect matrix  $\alpha $ is diagonal and set $ V = \delta \mu $ 
Under these conditions, it follows from (50) that 
\begin{align}
 \Gamma^T V = \Gamma^T \delta \mu = \delta \Gamma \delta^{-1}  \delta \mu = - \delta ( 1 + \alpha \theta ) \mu  = -  ( 1 + \alpha \theta ) \delta \mu = -( 1 + \alpha \theta ) V \nonumber
\end{align}
\begin{remark}
The condition imposed on $ \Gamma $  is not too restrictive. For example it is satisfied when $ \Gamma $ is symmetric or
when $\Gamma $ is a two by two matrix with off-diagonal elements having the same sign.
\end{remark}
Therefore, if we set $ R = -\Gamma^T V $ then R is a positive vector (determined up to a scalar multiplier) that satisfies condition (62). 
Finally, according to the constant effect demand model specification (13) we must have
\begin{eqnarray}
\log R \; \equiv \log \left (  p_0 I_0 \frac {\theta} {1 + \theta} \right ) \; = \; 
\alpha \log I_0 + ( 1 + \Gamma ) \log p_0  + \nonumber  \\ 
( \alpha \theta \mu  + ( 1 + \Gamma) \mu )\log \tau  \nonumber
\end{eqnarray}
The coefficient in front of $\tau $ vanishes by construction (50-51) and
we have additional conditions on initial prices and inventories   
\begin{align}
R = p_0 I_0 \frac{\theta}{1 + \theta}  \Leftrightarrow \log R =  \log p_0 + \log I_0 + \log \frac{\theta}{1 + \theta} \\
\log R = \alpha \log I_0 + ( 1 + \Gamma ) \log p_0 
\end{align}
Clearly, these conditions determine  $ p_0 $ completely and $ I_0 $ up to a scalar multiplier.   
%For example, this is the case when elasticity matrix $\Gamma $ is symmetric or if the number of items $ n = 2 $ and off-diagonal elements of $ \Gamma  $ are  of the same sign.  
%In any case, if $ V $ is a positive vector, then we can set  
% 
%\begin{equation}
% - V  \; = \;  \left( \Gamma^T \right)^{-1}  R  \; \Leftrightarrow \; -\Gamma^T V \; = \; R
% \; \Leftrightarrow  -\Gamma^T V  \; = \;  p_0 I_0 \frac {\theta} {1 + \theta} 
%\end{equation}
%for some appropriately chosen vector of prices $ p_0 $
%( and note, that according to the Remark 6   
%%and assuming that $ \alpha $ is a diagonal matrix with positive elements 
%we must have $ \Gamma^T V  < 0 $).
%Hence (43)  is an optimal solution of a markdown problem(5-6) as long as there exists a vector $ 0 < \mu < 1$ that satisfies (51).  
We thus have 
\begin{proposition} 
A solution for  markdown optimization problem (5-6) with constant symmetric elasticity matrix $ \Gamma, $ constant positive diagonal inventory effect matrix $ \alpha $ and boundary conditions determined by  (64-65) is given by price and inventory curves (43) if  
\begin{equation}
 [1_n - \alpha^{-1}(1 + \Gamma)] \frac{\theta}{1 + \theta} = 1_n 
\end{equation}
\end{proposition} 
% Proof. The sufficiency of condition (67) is already established. To prove the second statement of the theorem, suppose that (43) is a solution of (5-6). Then by 
% Theorem 2, a linear combination with constant coefficients of individual item revenues 
% \begin{equation}
%  f( t ) \; = \; \sum_{i}^n a_i R_i(t) 
%\end{equation}
%   computed along this solution must be constant. By (43) all the functions  
% $ R_i(t)  $ must be some non-negative powers of $ \tau  $ and from 
% again Theorem 2 it follows  that coefficients $ a_i $ are sums of rows of 
% $ \Gamma^{-1}   $ hence $ a_i < 0, $ for all $ i = 1,2,\cdots, ,n $ since $ \Gamma .$    
%  Collecting  non-constant terms in (67) we get
% \begin{equation}
%  f_1( t ) \; = \; \sum_{j} a_j \tau(t)^{\beta_j} = const ,  \; \beta_j > 0, \; j \in \{ 1,2,, \cdots ,n \}     
%\end{equation}
%Note that  $ f_1(T) = 0 $ (by the definition of $ \tau $)  and therefore $ f_1(t) $ must be identically 
% equal to $ 0 $.  This is in contradiction with $ a_j < 0 $ for al $ j = 1, \cdots , n.$ Hence we must assume that all individual revenues $ R_i $  are constant as is stipulated in (44). 
% The subsequent computations are clearly reversible.
\begin{remark}
The solution obtained in Proposition 2 in general can match only an arbitrary magnitude of the initial inventory vector, not its direction.  The reason for this drawback is excessive restrictiveness of condition (44). However, as we will see shortly the initial inventory condition can be satisfied if we are allowed to change the markdown time period.
\end{remark} 

\subsection{Concrete numerical example}
Here we present a straightforward numerical example mimicking the  Proposition 2 for a realistic two item markdown scenario. We set 
$$
\alpha = \left(
\begin{array}{cc}
0.5 & 0 \\
0 & 0.3 \\
\end{array}
\right) , \; \Gamma = \left(
\begin{array}{cc}
-2 & 0.25 \\
0.25 & -1.5 \\
\end{array}
\right),  \;  
I_0 = \left(
\begin{array}{cc}
200 \\
300 \\
\end{array}
\right)
$$
and compute vector $ \mu $ from equation (51)
$$
 \mu = \left(
\begin{array}{cc}
0.417 \\
0.505 \\
\end{array}
\right)
$$
Hence, 
$$
  \theta \mu  =  1 - \mu =
   \left(
\begin{array}{cc}
0.582 \\
0.494 \\
\end{array}
\right)   
$$
Ignoring seasonality effect we find optimal inventory curve
\begin{align}
    I(t) = I_0 \Bigl (1 - \frac{t}{T} \Bigr )^{\theta\mu}    
\end{align}
Since $ \alpha $ is diagonal and $ \Gamma $ is symmetric we find revenue rate to be
$$
R  = -\Gamma^T\mu = 
\left(
\begin{array}{cc}
0.708 \\
0.653 \\
\end{array}
\right) 
$$
Of course any scalar multiple  of the above vector will also satisfy relation (63) and what matters is the ratio $R_1 / R_2 = 1.084$. As was explained in  the previous section, the following conditions must be satisfied
\begin{align}
RT = p_0 I_0 \frac{\theta}{\theta + 1 }  \nonumber \\
\alpha \log I_0 + ( 1 + \Gamma ) \log p_0 = \log R  \nonumber
\end{align} 
The time period $T $ shows up here since our measure $ dt $ is not normalized 
as was the case for the revenue formula  (49).
These are two vector equations with unknown initial price  $ p_0, $
 magnitude of $ R $ and time period $ T .$ From here we find that 
\begin{align}
R(0) = 
\left(
\begin{array}{cc}
5.183 \\
4.781 \\
\end{array}
\right) , \;  p(0) = 
\left(
\begin{array}{cc}
3.424 \\
2.480 \\
\end{array}
\right) ,  \; 
T \approx 77
\end{align}  
with the ratio of the components of R being about $1.804$ as expected. 
It should be emphasized again (cf. Remark 10 above), that markdown solution presented in Proposition 2 is grossly overdetermined. A standard markdown problem specifies both markdown period and initial inventories.   
  
\section{In Conclusion: Direct Variational Methods}

We have shown that variational framework can be used  for qualitative analysis of retail pricing models.  It was also demonstrated that a "conservation law" described by Corollary 2 can be successfully applied to practical  revenue optimization problems.   However, analytical methodology can serve only as a guide to numerical analysis of the problems involved.
Closed form solutions similar to Propositions 1 and 2 
( cf. also \cite{Smith2010} and \cite{Vak} ) are very rare and can be obtained only in relatively simple cases.
Nevertheless, pricing models
with general constraints and boundary conditions can be  efficiently handled  in a general framework of direct 
 variational methods (see e.g. \cite{Boyd}, \cite{abs}, \cite{O}) while  qualitative results, such as Theorem 2 could serve as useful model specification guidelines.
Moreover, since direct variational methods essentially work with discrete data, the discrete nature of some business problems can be effectively handled by a variational framework as well.
A serious objection to the revenue optimization along these  lines is that  it is harder to estimate 
cross-elasticity matrix $ \Gamma $ than to compute an optimal price policy based on such
$ \Gamma .$  We hope, nevertheless, that presented results (e. g. formula (27)) could  provide some useful guidelines for (cross-)elasticity estimation.

%\begin{acknowledgements}
\paragraph{}
The authors are grateful to Andrew Vakhutinsky for stimulating discussions.
%\end{acknowledgements}


\begin{thebibliography}{}

\bibitem {Talluri} Kalyan T. Talluri, Garrettt J. Van Ryzin
\emph{ The theory and practice of Revenue Management, 2005} 

\bibitem{Smith1998} Stephen A. Smith , Dale D. Achabal \emph{ Clearance Pricing and Inventory Policies for
Retail Chains }, Management Science, Vol. 44, No. 3, pp. 285-300, March 1998  


\bibitem{Smith2010}  Stephen A. Smith  Clearance Pricing in Retail Chains, in 
Retail Supply Chain
Management, Quantitative Models and Empirical Studies, Springer, 271-291, 2009  

\bibitem{deim}  Klaus Deimling, Non Linear Functional Analysis,
Dover, 2010  

\bibitem{Ler} Lerner, A. P. ,  The Concept of Monopoly and the Measurement of Monopoly Power, The Review of Economic Studies,  Vol. 1, No. 3, 157-175, 1934

\bibitem{Vak} Andrew Vakhutinsky, Alex Kushkuley, Manish Gupte, Markdown Optimization under Inventory-Depletion Effect,   
Journal of Revenue and Pricing Management Vol. 11, 6, 632  - 644, 2012


\bibitem{Thijs}
Thijs ten Raa, Monopoly, Pareto and Ramsey Mark-ups,  Journal of Industry, Competition and Trade, 9:57–63, 57-63, 2009

\bibitem{abs}  Alexander Kushkuley, Su-Ming Wu, On Variational Approach to Multivariate Retail Pricing Models,  
Presentation at MSOM Conference, 2012
\bibitem{gelf} I. M. Gelfand and S. V. Fomin, Calculus of Variations, Dover, 1991 
% etc
\bibitem{horn} Roger A. Horn, Charles R. Johnson, Matrix Analysis, Cambridge University Press, 2013

\bibitem{Boyd}  John P. Boyd, \emph{Chebyshev and Fourier Spectral Methods},  Dover, 2000
\bibitem{O} O. von Strykand R. Bulirsch, Direct and indirect methods for trajectory  optimization, Annals of Operations Research 37 (1992), 357-373
\end{thebibliography}
\end{document}